\documentclass[11pt,leqno]{article}
\usepackage{amsfonts}
\usepackage[dvips]{graphicx}
\usepackage{latexsym}
\usepackage{amsmath}

\newcommand{\domek}{\rule{0.09in}{0.09in}}
\newcommand{\eps}{{\epsilon}_{M}}

\newenvironment{AMS}{\small\it AMS subject classification:}{ }
\newenvironment{keywords}{\small\it Keywords:}{ }
\newtheorem{theor}{Theorem}
\newtheorem{lemat}{Lemma}
\begin{document}

\author{{\small Alicja Smoktunowicz and Iwona Wr\'{o}bel}}
\title{{\Large On improving the accuracy of Horner's and Goertzel's
algorithms}}
\date{ {\small Faculty of Mathematics and Information Science, Warsaw
University of Technology, Pl.Politechniki 1, 00-661 Warsaw, Poland }\\
e-mail: smok@mini.pw.edu.pl, wrubelki@wp.pl}
\maketitle

\begin{abstract}
It is known that Goertzel's algorithm is much less numerically
accurate than the Fast Fourier Transform (FFT)(Cf. \cite{gen:69}).
In order to improve accuracy we propose modifications of both
Goertzel's and Horner's algorithms based on the divide-and-conquer
techniques. The proof of the numerical stability of these two
modified algorithms is given. The numerical tests in Matlab
demonstrate the computational advantages of the proposed
modifications. The appendix contains the proof of numerical
stability of Goertzel's algorithm of polynomial evaluation.
\end{abstract}

\noindent

\begin{AMS}
{\small 65F35, 65G50.}
\end{AMS}

\noindent


\section{Introduction}

The aim of this paper is to improve the accuracy of polynomial
evaluation, mainly Horner's and Goertzel's algorithms. Both,
Horner's and Goertzel's methods are frequently used in the
interpolation and approximation problems and in signal processing.
Goertzel's algorithm is implemented in Matlab, it's included in
the Signal Processing Toolbox. The function ''fft'' returns the
Discrete Fourier Transform (DFT) computed with a Fast Fourier
Transform (FFT) algorithm and the function ''goertzel'' computes
DFT of specific indices in a vector.

In this paper we consider more general case of evaluating a
polynomial
\begin{equation}
w(z)=\sum_{n=0}^{N}{a_n z^n},
\end{equation}
where $z\in \mathbb C$ and $a_0,\dots,a_N\in\mathbb C$.

Note that for $z=e^{i\xi}$ and $a_0,\dots,a_N,\xi\in\mathbb R$ we
have
\[
w(e^{i\xi})=\sum_{n=0}^N a_n \cos n\xi + i\sum_{n=0}^N a_n \sin
n\xi.
\]
DFT returns $y_k=w(z_k)$, $k=0,\dots,N$, where
$z_k=e^{i\xi_k}=e^{-\frac{2\pi i k}{N+1}}$ (Cf. \cite{loan:92}, p.
10).

It is observed (see ''help goertzel'' in Matlab Signal Processing
Toolbox) that compared with the Fast Fourier Transform algorithm
(FFT), Goertzel's algorithm is much less numerically accurate,
which can be visible especially for high-scale problems.


We propose the algorithm PEMA (Polynomial Evaluation Modified
Algorithm), which is based on the repetitive use of some algorithm
$W$ for evaluating polynomials. This algorithm can be e.g.
Horner's or Goertzel's scheme. The cost of PEMA is comparable to
the cost of $W$ and the error bound of PEMA may be significantly
smaller than the error bound of $W$. We prove that if $W$ is
stable then PEMA is also numerically stable (see section 3.2). In
practice, one should use only numerically stable algorithms.

We say that an algorithm of evaluating (1) is componentwise
backward stable with respect to the data $a_0,\dots,a_N\in\mathbb
C$ and $z \in\mathbb C$ if the value $\tilde w(z)$ computed by
this algorithm is an exact value of a polynomial for slightly
perturbed coefficients $a_n$ and $z$, i.e.
\begin{equation}
\tilde w(z)=\sum_{n=0}^{N}{\lbrack a_{n}(1+\mu_n)\rbrack\, \lbrack
z(1+\beta )\rbrack^n}, \hskip0.4cm \vert \mu_n \vert \leq A_N
\eps, \hskip0.4cm \vert \beta \vert \leq Z_N \eps,
\end{equation}

\noindent where $A_N$ and $Z_N$ are modestly growing functions of
$N$ and $\eps$ is the machine precision.

Throughout the paper we assume that the coefficients of a
polynomial $w(z)$ are complex.

In the error analysis of PEMA we consider perturbations not only
of polynomial coefficients, but also of $z$. Notice that usually
the exact value of $z$ is not known, e.g. $z$ is given as
$z=e^{i\xi}$. Then $z=c+i\,s$, $c=\cos\xi$, $s=\sin\xi$ and
$\tilde c =c+\Delta c$, $\tilde s =s+\Delta s,\,$ $\vert \Delta c
\vert, \: \vert \Delta s \vert \leq \nu \,\eps$, where $\nu$ is
small. Then the perturbed value $\tilde z$ can be written as
$\tilde z= z\,(1+\eta)$, $\vert \eta \vert \leq \sqrt 2 \,\nu
\,\eps$.

Then with help of Taylor expansion (2) leads to
\[
\ \mid \!\tilde w(z) - w(z)\! \mid \leq \mid \!\beta\! \mid \, \mid\! z\, w^\prime(z)\! \mid + \sum_{n=0}^{N}\mid \!a_{n}\!\mid \,\mid\!\delta_n
\!\mid \,\mid \!z \!\mid^n + \,\mathcal O(\eps^{\,2}),
\]

\noindent and further
$$\ \mid \!\tilde w(z) - w(z)\! \mid \leq \eps (A_N \sum_{n=0}^{N}\mid \!a_{n}\!\mid \,\mid \!z
\!\mid^n + \,Z_N\! \mid\! z \! \mid\,\mid\! w^\prime(z)\! \mid)
+\,\mathcal O(\eps^{\,2}). $$

Numerical stability  of Horner's algorithm was first given by
Wilkinson (Cf. \cite{wilk:63}, pp. 36-37, 49-50) who proved that
$Z_N=0$ and $A_N \approx 2 N$, provided that the data $a_0, \dots,
a_N \in \mathbb R$ and $z \in \mathbb R$ are exactly representable
in floating point arithmetic (fl). Despite of a bad reputation of
Goertzel's algorithm as a method of computing Fourier series
$\sum_{n=0}^{N}{a_n \cos n\xi}$ and $\sum_{n=1}^{N}{a_n \sin
n\xi}$ with respect to the data $a_0, \dots, a_N \in \mathbb R$
and a given argument $\xi \in \mathbb R$ (Cf. \cite{stoer:80}, pp.
84-88, \cite{gen:69}, \cite{newbery:73}, \cite{oliver:77}) we
prove that Goertzel's algorithm is numerically stable in a sense
(2). The respective constants are $Z_N=0$ and $A_N$ is of order
$N^2$, provided that the data $a_0, \dots, a_N$ and $z$ are
exactly representable in fl (see Theorem 2 and Table 0).

In order to improve accuracy we propose modifications of both
Goertzel's and Horner's algorithms based on the divide-and-conquer
techniques. The idea is not quite new, there are numerous
divide-and-conquer parallel algorithms for polynomial evaluation
(Cf. \cite{gol:93}, p. 70). The goal of our work is to split a
polynomial in ''the proper way'' in order to refine results. We
show that the constants $A_N$ and $Z_N$ in (2) can be
significantly decreased, in comparison with the classical Horner's
and Goertzel's algorithms, which is of great importance for large
$N$ (see Table 0 in section 3), e.g. for $N=2^p$ our
divide-and-conquer algorithm PEMA results in $A_N$ of order $log_2
N$ and $Z_N$ of order unity.

Tests included in section 4 confirm theoretical results. We also
implemented Reinsch's modification of Goertzel's algorithm (Cf.
\cite{stoer:80}, pp. 86-88) for evaluation of (1), but it turned
out that the numerical results were comparable to these given by
standard Goertzel's algorithm. For this reason we don't include
them in section 4 devoted to numerical experiments.

\section{Classical polynomial evaluation schemes}

\noindent The Horner scheme is the standard method for evaluation
of a polynomial (1)
at a given point $z \in \mathbb C$. We assume that
$a_0,\dots,a_N\in\mathbb C$. We write $w(z)$ as follows
$$\ w(z)=a_0+ z (a_1 + z(\ldots + z(a_{N-1}+z a_N)\ldots )).$$

\bigskip
\noindent {\large \bf Algorithm 1 (Horner's rule)}

\begin{flushleft}
$w:=0$ \\
for \hskip.2cm $n=N, N-1, \ldots, 0 $ \\
\hskip.8cm $w:=a_{n}+ z \, w$ \\
end \\
$w(z):=w$
\end{flushleft}

The complexity of Horner's algorithm $C_N(H)$, counted as a number
of multiplications is equal to $N$, which gives in general
$C_N(H)=4N$ real multiplications. We assume that the product of
two complex numbers is computed in a natural way and in
consequence one complex multiplication is equivalent to four real
ones.

The idea of Goertzel's algorithm is different. Suppose $z=x+i y$.
Divide a polynomial $w(\lambda)=\sum\limits_{n=0}^{N}\,\, a_n
\lambda^n$ by a quadratic polynomial $(\lambda- z)
(\lambda-\overline{z})= \lambda^2 - \hat p \lambda -\hat q$ with
real coefficients $\hat p$ and $\hat q$, where $\hat p=2 x$ and
$\hat q= - {\vert z \vert}^2$. Then
\[
w(\lambda)=(\lambda- z) (\lambda-\overline{z}) \sum_{n=2}^{N}{\
b_{n} \lambda^{n-2}} + b_0+ b_1\,\lambda
\]
\noindent and, consequently,
$w(z)= b_0+ b_1\,z.$ \noindent This leads to the following

\bigskip
\noindent {\large \bf Algorithm 2 (Goertzel's algorithm)}
\nopagebreak
\begin{flushleft}
$\hat p:= 2 x$\\
$\hat q:=- (x^2+y^2)$\\
$b_{N+1}:=0$\\
$b_{N}:=a_{N}$\\
for \hskip.2cm $n=N-1, \ldots, 1 $\\
\hskip.8cm $b_{n}:=a_{n}+ \hat p \,\, b_{n+1} + \hat q \,\, b_{n+2} $ \\
end \\
$u:= (a_{0}+ x \, b_1 + \hat q\, b_2)$ \\
$v:= y \, b_1 $ \\
$w( z ):= u + i \,v $ \\
\end{flushleft}

In general, the number of real multiplications needed by
Goertzel's method is the same as those needed by Horner's
algorithm. However, in special cases each of these algorithms can
be less expensive than the other. For example, for $z \in \mathbb
R$ Goertzel's algorithm is twice as expensive as Horner's rule
regardless of the polynomial coefficients. On the other hand
consider the case of polynomial with real coefficients and $z \in
\mathbb{C}$ , $\vert z \vert = 1$. Then all $b_n$ are real, $\hat
q=1$ and
$\ b_{n}=a_{n}+ \hat p \,\, b_{n+1} + \,\, b_{n+2}$ for $n=N-1,
\dots, 1$.

\noindent The complexity of Goertzel's method reduces to $N$ while
the cost of Horner's rule is still $4N$.

Note that if $z=1$, then $w(z)=\sum_{n=0}^{N}{a_{n}}$ and Horner's
rule is nothing else but a backward summation.

We now derive an algorithm based on the divide-and-conquer
technique.

\section{A new polynomial evaluation modified algorithm (PEMA)}

Suppose a polynomial $w(z)$ is given by
$w(z)=\sum\limits_{n=0}^{N}\,\, a_nz^n$ where $N=s^p$ and $s>1$.
We can write $w(z)$ in the following form:
\begin{multline}
w(z)=\{a_0+a_1z+\dots+a_{s-1}z^{s-1}\}+\{a_s+a_{s+1}z+\dots+a_{2s-1}z^{s-1}\}z^s+\dots+\nonumber\\
+\{a_{(s^{p-1}-1)s}+a_{(s^{p-1}-1)s+1}z+\dots+a_{s^p-1}z^{s-1}\}(z^s)^{s^{p-1}-1}+a_{s^p}(z^s)^{s^{p-1}}=\nonumber\\
=a_0^{(1)}+a_1^{(1)}z_1+a_2^{(1)}z_1^2+\dots+a_{s^{p-1}}^{(1)}z_1^{s^{p-1}}=\sum_{j=0}^{s^{p-1}}a_j^{(1)}z_1^j,\nonumber
\end{multline}
where $z_1=z^s$, $a_j^{(1)}=\sum_{k=0}^{s-1}a_{js+k}^{(0)}z^k$,
$j=0,1,\dots,s^{p-1}-1$, $a_{s^{p-1}}^{(1)}=a_{s^p}^{(0)}$, and
$a_j^{(0)}=a_j$ for $j=0,1,\dots, N$.

Now we can interpret $\sum_{j=0}^{s^{p-1}}a_j^{(1)}z_1^j$ as a
polynomial of variable $z_1$ with the coefficients $a_j^{(1)}$ and
proceed in the same manner as before. We continue this process and
for $m=0,1, \ldots, p-1$ write $w(z)$ as follows
$$\ w(z)= \sum\limits_{j=0}^{s^{p-m}}\,\, a_j^{(m)}z_m^j, $$

\noindent where $z_0=z$ and $z_m=z_{m-1}^s$ for $m=1,2,\dots,p-1$.

\medskip

It is easy to prove that for $m=1,2,\dots,p-1$ and
$j=0,1,\dots,s^{p-m}\!-1$
\begin{equation}
a_j^{(m)}=\sum\limits_{r=0}^{s^m-1}\,\,a_{js^m+r}\,z^r.
\end{equation}
For complexity and computational accuracy reasons we don't
evaluate (3) directly, by Horner or Goertzel algorithm for
polynomial of variable $z$ and degree $s^m-1$, but use the
relation
\[
a_j^{(m)}=\nolinebreak\sum\limits_{k=0}^{s-1}\,\,a_{js+k}^{(m-1)}z_{m-1}^k.
\]
Notice that $a_j^{(m)}$ is a polynomial of variable $z_{m-1}$ and
degree $s-1$.

\smallskip
More precisely, given an algorithm $W$ for evaluating polynomials,
e.g. Horner's or Goertzel's algorithm, we produce a new
divide-and-conquer algorithm.

\bigskip

\noindent {\large \bf Algorithm 3 (PEMA)}
\medskip

This algorithm uses the divide-and-conquer method to compute
$w(z)$ where $z \in \mathbb R$ or $z \in \mathbb C$. The
coefficients $a_n$ may be either complex or real.

\begin{enumerate}
\item
\begin{flushleft}
$ z_0=z $ \\
\smallskip
$ a_{j}^{(0)}=a_j $ \hskip.2cm for \hskip.2cm $j=0,1,\dots,N$
\end{flushleft}
\item
\begin{flushleft}
for \hskip.2cm $m=1, \dots ,p-1 $ \\
\smallskip
\hskip.5cm $z_m=z_{m-1}^s $ \\
\smallskip
\hskip.5cm $a_{s^{p-m}}^{(m)}=a_{s^{p-(m-1)}}^{(m-1)}=a_N $ \\
\smallskip
\hskip.5cm for \hskip.2cm $j=0,1,\dots,s^{p-m}-1 $ \\
\hskip1cm compute \hskip.2cm $a_j^{(m)}=\nolinebreak\sum\limits_{k=0}^{s-1}\,\,a_{js+k}^{(m-1)}z_{m-1}^k $ \hskip.2cm by algorithm $W$ \\
\hskip.5cm end \\
\smallskip
end
\end{flushleft}
\item
compute \hskip.2cm $w(z)=\sum\limits_{j=0}^s\,\,a_j^{(p-1)}z_{p-1}^j $ \hskip.2cm by algorithm $W$
\end{enumerate}

\noindent Note that $p=1$ implies $N=s$ and PEMA is nothing else
but $W$ applied to $w(z)$.

PEMA is an extension of a summation algorithm proposed in
\cite{jank:83}. For $N=2^p$ and $z=1$ PEMA coincides with the
log-sum algorithm.

\vskip.5cm

\subsection{Total cost of PEMA}

\medskip
\nopagebreak \noindent Suppose the complexity of the algorithm $W$
is $C_N=bN$, $b=const$, i.e. $W$ needs $C_N$ multiplications to
compute $w(z)= \sum\limits_{n=0}^{N}\,\, a_nz^n$. We give a
formula for complexity of PEMA valid under assumption that $z_m$
is computed in a natural way, (see section 3.2):
$$\ C(PEM\!A)\!=\!\!\!\sum\limits_{m=1}^{p-1}\lbrace\sum\limits_{j=0}^{s^{p-m}\!-1}\!\!C_{s-1}+(s-1)\rbrace+C_s\!=\!(s-1)(p-1)+C_s+C_{s-1}\frac{s(s^{p-1}\!-\!1)}{s-1}. $$

\noindent According to this formula the complexity of PEMA with
Horner is equal to $C_N+(s-1)(p-1)$. Very often the latter term is
not significant in comparison with $C_N$.
\medskip

\noindent {\sc \bf Remark.} Each $a_j^{(m)}\!$,
$j=0,1,\dots,s^{p-m}-1$ can be computed independently. It's a big
advantage of PEMA because of possibility of parallel
implementation, which can be useful especially for really large
problems.

\subsection{Error analysis of PEMA}

\noindent We consider complex arithmetic (cfl) implemented using
standard real arithmetic with machine precision $\eps$. Then
\begin{equation}
cfl(x + y) = (x + y) (1+\delta), \hskip.6cm \mid \delta \mid \leq
\eps \hskip.6cm \mathrm{for} \;\;x,y \in \mathbb C
\end{equation}

\noindent and provided that the product $xy$ is computed using an
ordinary algorithm we have (Cf. \cite{krysia:76})
\begin{equation}
cfl(xy) = (xy) (1+\eta), \hskip.6cm \mid \eta \mid \leq c\, \eps,
\end{equation}

\noindent where
\begin{equation}
c=\left\lbrace\begin{array}{cl} 1 & \mathrm{for} \hskip.5cm x,y
\in \mathbb R \hskip.4cm \mathrm{or} \hskip.3cm x\in \mathbb R,
\hskip.15cm y \in \mathbb C
\\
1+\sqrt 2 & \mathrm{for} \hskip.5cm x,y \in \mathbb C.
\end{array}\right.
\end{equation}

\noindent The value $z_m=z_{m-1}^s$ is determined in a natural way
by computing the consecutive powers of $z_m$, i.e. $z_{m-1}$,
$z_{m-1}^2$, $\dots$, $z_{m-1}^s$.


\noindent Then
\begin{equation}
\tilde z_m = cfl(\tilde z_{m-1}^s) = \tilde z_{m-1}^s
(1+\delta_m), \hskip.6cm \mid \delta_m \mid \leq (s-1)\,c\, \eps
+\mathcal O(\eps^{\,2}).
\end{equation}



Now we are in a position to give the error analysis of the PEMA
algorithm. For simplicity we assume that $a_0,\dots,a_N$ and $z$
are represented exactly in cfl and that $s$ and $p$ are fixed,
$N=s^p$. We also assume that the result given by the algorithm $W$
of evaluating $w(z)=\sum_{n=0}^{N}{a_n z^n}$ in cfl satisfies
\begin{equation}
\tilde w(z)=\sum\limits_{n=0}^N\,\,a_n(1+\Delta_n)\,z^n,
\hskip.7cm \mid \Delta_n \mid \leq A_N \, \eps,
\end{equation}

\noindent where $A_N$ is an increasing function of $N$. $W$ in
PEMA can be Horner's or Goertzel's rule. For detailed information
on $A_N$ see (37).

For $m=1,\dots,p-1$ and $j=0,1,\dots,s^{p-m}-1$ the values $\tilde
a_j^{(m)}$, computed in cfl, can be written as follows
\begin{equation}
\tilde a_j^{(m)}=\sum\limits_{k=0}^{s-1}\,\,\tilde
a_{js+k}^{(m-1)}(1+\Delta_{j,k}^{(m)})\,\tilde z_{m-1}^k,
\hskip.7cm \mid \Delta_{j,k}^{(m)} \mid \leq A_{s-1} \, \eps.
\end{equation}

\noindent The formula (7) allows us to write $\tilde z_{m-1}$ in
the following way
\begin{equation}
\tilde z_{m-1} = \lbrack z(1+\gamma_m)\rbrack^{s^{m-1}},
\hskip.7cm
1+\gamma_m=\prod\limits_{t=1}^{m-1}(1+\delta_t)^\frac{1}{s^t}.
\end{equation}
\noindent From (7) we obtain an upper bound for $\mid \gamma_m
\mid$
\[
\mid \gamma_m \mid \leq (s-1)\, c\, \eps \,
\sum_{t=1}^{m-1}{\frac{1}{s^t}}\,\, +\mathcal O(\eps^{\,2}).
\]
Thus
\begin{equation}
\mid \gamma_m \mid \leq c\, \eps +\mathcal O(\eps^{\,2}).
\end{equation}

\begin{lemat}Assume that $cfl(z)=z$ and $cfl(a_n)=a_n$, $n=0,\dots,N$
and $N=s^p$. Suppose that $A_N \,\eps \leq 0.1$ and that {\rm
(7-9)} hold. Then for $m=1,\dots, p-1$ and $j=0,1,\dots,s^{p-m}-1$
\begin{equation}
\tilde a_j^{(m)}=\sum\limits_{r=0}^{s^m-1}\,\,\lbrack
a_{js^m+r}(1+\eta_{j,r}^{(m)})\rbrack \, \lbrack
z(1+\gamma_m)\rbrack^r
\end{equation}
\noindent where
\begin{equation}
\mid \gamma_m \mid \leq c\, \eps +\mathcal O(\eps^{\,2}),
\end{equation}
\noindent and
\begin{equation}
\mid \eta_{j,r}^{(m)}\mid \leq m\,d \, \eps +\mathcal
O(\eps^{\,2}), \hskip.7cm d= A_{s-1}+(s-1)\,c,
\end{equation}
\noindent where $c$ is defined by {\rm(6)}.
\end{lemat}

\noindent {\it Proof.} Let $m=1$. Then from (9) it follows that
\begin{equation}
\tilde a_j^{(1)}=\sum\limits_{k=0}^{s-1}
a_{js+k}(1+\Delta_{j,k}^{(1)})\,z^k, \hskip.7cm \mid
\Delta_{j,k}^{(1)} \mid \leq A_{s-1} \, \eps,
\end{equation}
\noindent which can be rewritten in the following form
$$\ \tilde a_j^{(1)}=\sum\limits_{k=0}^{s-1}
a_{js+k}\frac{(1+\Delta_{j,k}^{(1)})}{(1+\gamma_1)^k}\,\lbrack
z(1+\gamma_1)\rbrack^k = \sum\limits_{k=0}^{s-1} \lbrack
a_{js+k}(1+\eta_{j,k}^{(1)})\,\lbrack z(1+\gamma_1)\rbrack^k,
$$
\noindent where
$$\ 1+\eta_{j,k}^{(1)} = \frac{1+\Delta_{j,k}^{(1)}}{(1+\gamma_1)^k}. $$
\noindent From this we obtain
$$\ \mid \eta_{j,k}^{(1)}\mid
\leq \mid \Delta_{j,k}^{(1)} \mid + \,k\mid \gamma_1 \mid + \,
\mathcal O(\eps^{\,2}). $$

\noindent Now using (11), (15), the definition of $d$ in (14) and
the fact that $k\leq s-1$ we have
$$\ \mid
\eta_{j,k}^{(1)}\mid \leq (A_{s-1}+(s-1)\,c)\, \eps +\mathcal
O(\eps^{\,2}) = d\, \eps +\mathcal O(\eps^{\,2}). $$

\noindent In the same manner, using  the equality
\[
[ z(1+\gamma_m) ] = [ z(1+\gamma_{m-1}) ]\,
(1+\delta_{m-1})^{\frac{1}{s^{m-1}}}
\]
we get
$\mid \eta_{j,k}^{(m)}\mid \leq m\,d\,
\eps +\mathcal O(\eps^{\,2})$, which is the desired
conclusion. \domek

\begin{theor} Under the assumptions of Lemma 1 the value $\tilde w(z)$ computed by PEMA
satisfies
$$\ \tilde w(z)= \sum\limits_{n=0}^N\,\,\lbrack a_n(1+\Delta_n)\rbrack\,\lbrack
z(1+\beta)\rbrack ^n,$$
\noindent where
$$\ \mid \beta \mid = \mid \gamma_p \mid
\leq c\, \eps +\mathcal O(\eps^{\,2}), \hskip.6cm \mid \Delta_n
\mid \leq p\,(A_s+s\,c)\, \eps +\mathcal O(\eps^{\,2}).$$
\end{theor}

\noindent {\it Proof.} Let $m=p-1$. Lemma 1 yields
$$\ \tilde a_j^{(p-1)}=\sum\limits_{r=0}^{s^{p-1}-1}
\lbrack a_{js^{p-1}+r}(1+\eta_{j,r}^{(p-1)})\rbrack\,\lbrack
z(1+\gamma_{p-1})\rbrack^r, $$

\noindent where
$$\ \mid \eta_{j,r}^{(p-1)}\mid \leq (p-1)\,d\, \eps
+\mathcal O(\eps^{\,2}).
$$

\noindent By assumptions, we have

$$\ \tilde w(z)= \sum\limits_{j=0}^s\,\,\tilde
a_j^{(p-1)}(1+\xi_{j}^{(p)})\,\tilde z_{p-1}^j, \hskip.6cm \mid
\xi_j^{(p)} \mid \leq A_s\, \eps.$$

\noindent This gives immediately the assertion of the theorem.
\domek

So, if the algorithm $W$ satisfies (8), PEMA is numerically stable
in a sense (2).

\medskip
\begin{center}
\noindent {\footnotesize Table 0: Constants $A_N$ and $Z_N$ for
all algorithms.}
\end{center}
\nopagebreak
\begin{center}
\begin{tabular}{|c|c|c|}
\hline Algorithm & $A_N$ & $Z_N$ \\
\hline Horner & $(c+1)N$ & 0 \\
\hline Goertzel & $10N^2$ & 0 \\
\hline PEMA(Horner) & $ps(2c+1)$ & $c$\\
\hline PEMA(Goertzel) & $10ps^2+psc$ & $c$\\
\hline
\end{tabular}
\end{center}

\noindent $N$ is the degree of the polynomial, $c=1$ for real
coefficients $a_0,\dots,a_N$ and $c=1+\sqrt 2$ for complex $a_n$.
Here $p$ and $s$ are the parameters of PEMA, $N=s^p$. Note that
for $N=2^p$, the partial polynomials in PEMA are of degree $1$ and
$A_N$ is of order $\log_2N$, which is a significant improvement
when compared with the standard versions of both algorithms.

\section{Numerical tests}

\noindent This paragraph contains the results of the tests
performed in Matlab, version 6.1.0450 (R12.1) with machine
precision $\eps \approx 2.2 \cdot 10^{-16}$. We implemented all
methods and compared the results they gave. Of course, it would be
the most natural to compare the result given by each of the
methods with the exact one. However, there are obvious
obstructions, i.e. for fractional polynomial coefficients or the
point $z$ there is no way to obtain the exact value of $w(z)$. To
deal with these difficulties we used the Matlab function ''fft'',
which is perfectly stable (for details see \cite{loan:92}, pp.
22-45). The function $yf\!ft=f\!ft(a)$ computes Fourier
coefficients, namely $y_k=w(z_k)$, $k=0,\dots,N$, where $w(z)$ is
the polynomial (1) and $z_k=\omega^k$, $\omega = e^{-\frac{2\pi
i}{N+1}}$ is the $(N+1)$st root of unity: $\omega ^{N+1}=1$. The
values $z_k$ were computed by the Direct Call algorithm, i.e.
$z_k:=\cos{(kt)}-i\sin{(kt)}$, $t=\frac{2\pi}{N+1}$, which is
known to be very accurate (Cf. \cite{loan:92}, pp. 23-24).

We computed the relative error
\begin{equation}
error=\frac{\Vert y-yf\!ft\Vert_2}{\Vert yf\!ft\Vert_2}
\end{equation}
\noindent where $y$ denotes the vector of results given by
Horner's, Goertzel's or PEMA algorithm for a certain set of points
$\{z_j\}$ and $yf\!ft$ is the result given by the function ''fft''
for the same set of points, namely for $z_j$, where $j\in
\{0,1,9,99,199,256,299,399,499,699\}$. The parameter $p$ in PEMA
(see section $3$) was equal to $2$, namely $N=s^2$ (i.e. $s=\sqrt
N$).


The function ''fft'' can be used provided that $\vert z\vert \le
1$. In general this condition is not needed, all algorithms,
namely Goertzel's, Horner's and both versions of PEMA work for any
$z\in \mathbb C$.

Figure $1$ describes the results for Goertzel's algorithm and PEMA
with Goertzel's method applied to polynomials with random
coefficients.

Both graphs illustrate the logarithm of $error$ (16) plotted
against the logarithm of the polynomial degree $n$, which varies
between $2^{10}$ and $2^{22}$. The lower graph represents results
given by PEMA, while the upper one these given by the standard
Goertzel's algorithm.

\vspace*{0.5cm}

\noindent \hspace*{0.85cm}
\includegraphics[width=10cm]{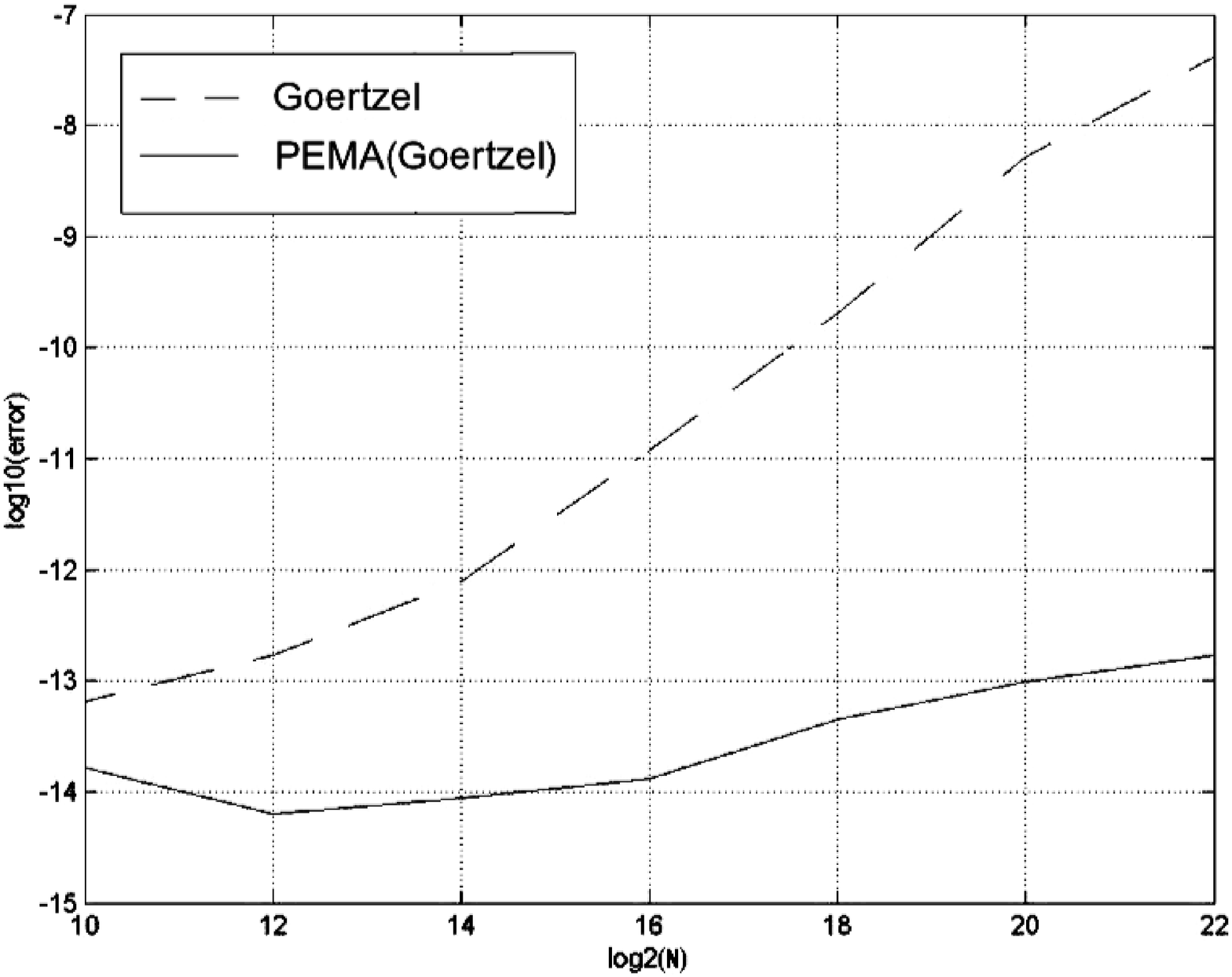}
\nopagebreak

\noindent \hspace*{1.8cm} {\footnotesize Figure 1: Relative errors of Goertzel's and PEMA algorithms} \\
\hspace*{3.25cm} {\footnotesize for polynomials with random
coefficients.}

\bigskip

Figure $2$ describes similar results for a family of polynomials
with coefficients given by the formula $a_{k}=f(t_{k})$ where
$t_k=0.001k$, $k=0,\dots, N$ and $f(t)=sint+sin100t+sin1000t$. As
before the lower and the upper graphs represent results given by
PEMA and the standard Goertzel's algorithm, respectively.

\vspace*{0.5cm}
\noindent \hspace*{0.85cm}
\includegraphics[width=10cm]{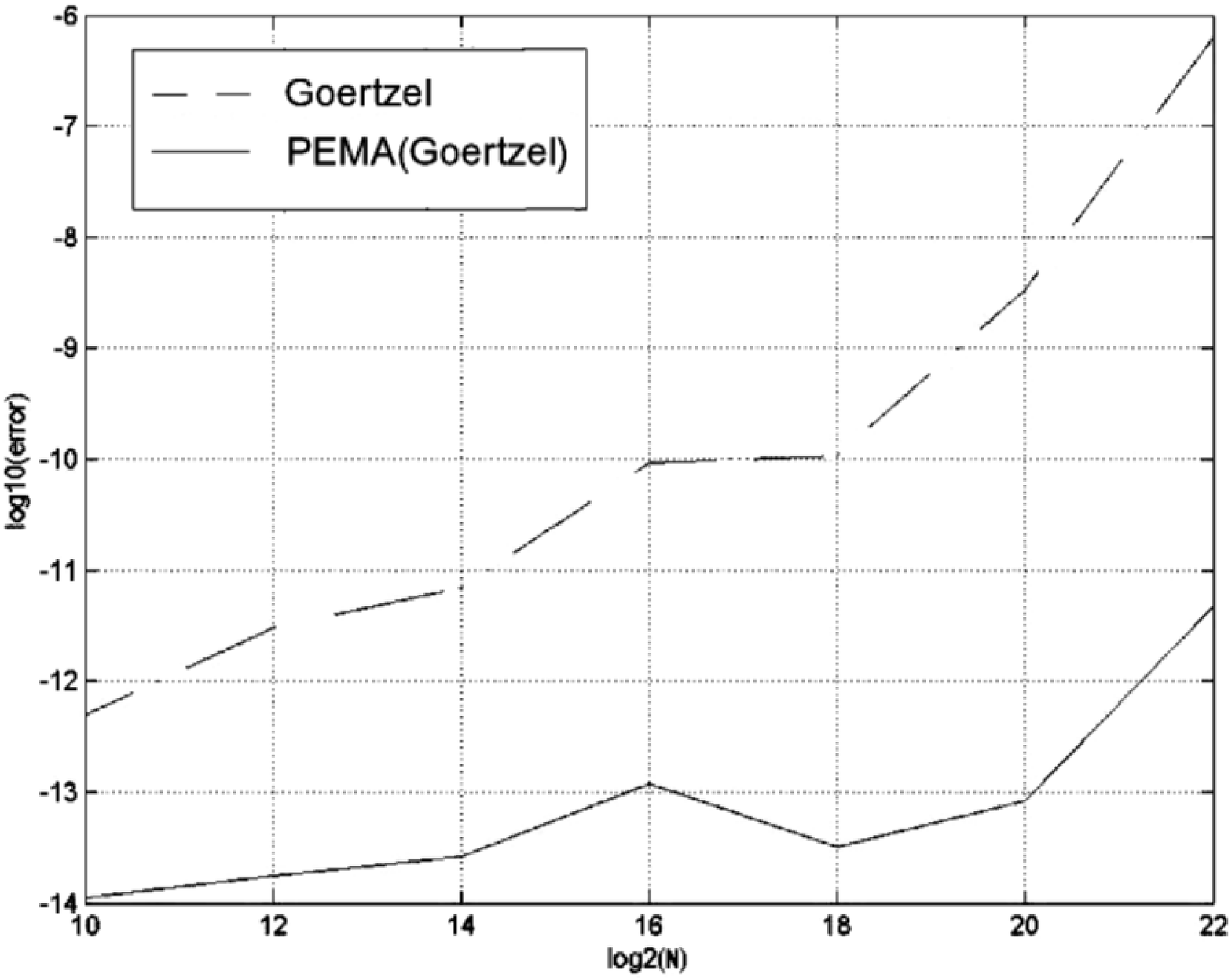}
\nopagebreak

\noindent \hspace*{1.8cm} {\footnotesize Figure 2: Relative errors of Goertzel's and PEMA algorithms} \\
\hspace*{3.25cm} {\footnotesize for polynomials with coefficients $a_{k}=f(t_{k})$} \\
\hspace*{3.25cm} {\footnotesize where
$f(t)=sint+sin100t+sin1000t$.}

\bigskip
\medskip

Figure 3 illustrates analogous results for polynomials with
coefficients $a_{k}=\sqrt{k}$. And again the lower graph
represents results given by PEMA.

\vspace*{0.5cm}
\noindent \hspace*{0.85cm}
\includegraphics[width=10cm]{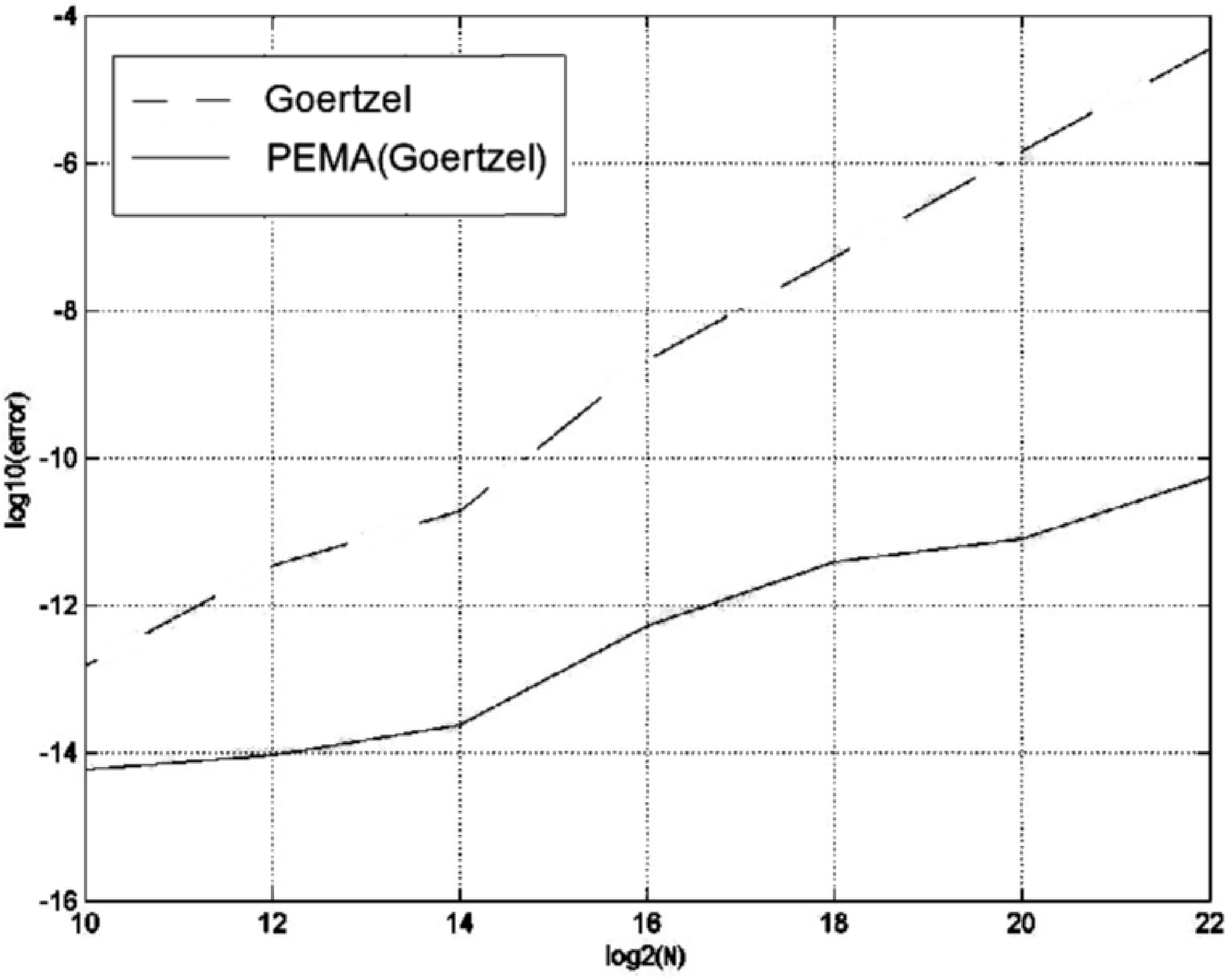}
\nopagebreak

\noindent \hspace*{1.8cm} {\footnotesize Figure 3: Relative errors of Goertzel's and PEMA algorithms} \\
\hspace*{3.25cm} {\footnotesize for polynomials with coefficients $a_{k}=\sqrt{k}$}. \\

\medskip

Tables $1-3$ contain values of $error$ (16) for each method and
for polynomials of coefficients given in description above each
table. $N$ is the polynomial degree. The second and the third
columns contain results given by Horner's rule and the version of
PEMA algorithm with Horner's rule, respectively. Data in the last
two columns is results given by Goertzel's algorithm and PEMA with
Goertzel's algorithm. This data was used to create figures $1-3$.

Note that although Goertzel's algorithm gives large errors for
large N, PEMA using Goertzel's algorithm has much smaller errors;
they are comparable with the errors obtained using Horner's
algorithm, or PEMA with Horner's algorithm.

\bigskip
\smallskip
\begin{center}
\noindent {\footnotesize Table 1: Relative errors of Goertzel's, Horner's and both versions of PEMA algorithms} \\
\hspace*{-4cm} {\footnotesize for polynomials with random
coefficients.}
\end{center}
\medskip
\begin{center}
\begin{tabular}{|c|c|c|c|c|}
\hline $N$ & Horner & $\!\!$PEMA(Horner)$\!\!$ & Goertzel &
$\!\!$PEMA(Goertzel)$\!\!$ \\ \hline $2^{10}$ & $1.6396e-014$ &
$1.6597e-014$ & $6.4827e-014$ & $1.6614e-014$
\\ \hline
$2^{12}$ & $6.4839e-015$ & $6.2312e-015$ & $1.7241e-013$ &
$6.3318e-015$
\\ \hline
$2^{14}$ & $6.4597e-015$ & $8.8147e-015$ & $7.8870e-013$ &
$8.8450e-015$
\\ \hline
$2^{16}$ & $1.0575e-014$ & $1.2730e-014$ & $1.1884e-011$ &
$1.3035e-014$
\\ \hline
$2^{18}$ & $3.0060e-014$ & $4.3985e-014$ & $2.0332e-010$ &
$4.4917e-014$
\\ \hline
$2^{20}$ & $7.1352e-014$ & $7.6212e-014$ & $5.2591e-009$ &
$9.7373e-014$
\\ \hline
$2^{22}$ & $1.1814e-013$ & $1.5060e-013$ & $4.1586e-008$ &
$1.7229e-013$ \\ \hline
\end{tabular}
\end{center}
\bigskip
\bigskip
\begin{center}
\noindent {\footnotesize Table 2: Relative errors of Goertzel's, Horner's and both versions of PEMA algorithms} \\
\hspace*{-3.8cm} {\footnotesize for polynomials with
coefficients $a_{k}=f(t_{k})$}\\
\hspace*{-4.3cm} {\footnotesize where
$f(t)=sint+sin100t+sin1000t$.}
\end{center}
\medskip
\noindent
\begin{center}
\begin{tabular}{|c|c|c|c|c|}
\hline $N$ & Horner & $\!\!$PEMA(Horner)$\!\!$ & Goertzel &
$\!\!$PEMA(Goertzel)$\!\!$ \\ \hline $2^{10}$ & $2.1321e-015$ &
$1.0999e-014$ & $4.9016e-013$ & $1.1313e-014$
\\ \hline
$2^{12}$ & $4.3372e-015$ & $1.5549e-014$ & $3.0108e-012$ &
$1.7462e-014$
\\ \hline
$2^{14}$ & $9.7481e-015$ & $2.5365e-014$ & $7.0483e-012$ &
$2.6262e-014$
\\ \hline
$2^{16}$ & $3.2760e-014$ & $1.0139e-013$ & $9.2595e-011$ &
$1.1973e-013$
\\ \hline
$2^{18}$ & $1.6703e-014$ & $1.6408e-014$ & $1.0737e-010$ &
$3.2052e-014$
\\ \hline
$2^{20}$ & $8.1798e-014$ & $6.0448e-014$ & $3.3356e-009$ &
$8.3002e-014$
\\ \hline
$2^{22}$ & $2.6576e-011$ & $3.9179e-011$ & $6.4574e-006$ &
$4.8041e-011$ \\ \hline
\end{tabular}
\end{center}

\bigskip
\bigskip
\begin{center}
\noindent {\footnotesize Table 3: Relative errors of Goertzel's, Horner's and both versions of PEMA algorithms} \\
\hspace*{-4cm} {\footnotesize for polynomials with coefficients
$a_{k}=\sqrt{k}.$}
\end{center}
\medskip
\noindent
\begin{center}
\begin{tabular}{|c|c|c|c|c|}
\hline $N$ & Horner & $\!\!$PEMA(Horner)$\!\!$ & Goertzel &
$\!\!$PEMA(Goertzel)$\!\!$ \\ \hline $2^{10}$ & $5.6281e-015$ &
$5.6566e-015$ & $1.5038e-013$ & $5.9073e-015$
\\ \hline
$2^{12}$ & $8.0767e-015$ & $8.1583e-015$ & $3.4601e-012$ &
$9.3555e-015$
\\ \hline
$2^{14}$ & $1.8735e-014$ & $1.8795e-014$ & $1.9228e-011$ &
$2.3707e-014$
\\ \hline
$2^{16}$ & $1.7620e-013$ & $4.7930e-013$ & $2.1008e-009$ &
$5.2504e-013$
\\ \hline
$2^{18}$ & $1.1682e-012$ & $3.5980e-012$ & $5.1238e-008$ &
$3.8532e-012$
\\ \hline
$2^{20}$ & $8.4972e-012$ & $6.1673e-012$ & $1.4374e-006$ &
$8.1276e-012$
\\ \hline
$2^{22}$ & $5.1749e-011$ & $4.1890e-011$ & $3.5824e-005$ &
$5.3874e-011$ \\ \hline
\end{tabular}
\end{center}

\bigskip
\bigskip
{\large \bf Appendix. Error analysis of Goertzel's algorithm}
\bigskip

\noindent Now we turn our attention to numerical analysis of
Goertzel's algorithm. Goertzel's method is a special case of
Clenshaw's algorithm (Cf. \cite{clen:55}, \cite{gen:69},
\cite{smokt:02}). Our results are similar in spirit to these given
by Gentleman \cite{gen:69}, who gave a floating-point error
analysis of Goertzel's algorithm for computing Fourier
coefficients $\sum_{n=0}^{N}{a_n \cos n\xi}$ and
$\sum_{n=1}^{N}{a_n \sin n\xi}$ with respect to the data $a_0,
\dots, a_N$ and a given argument $\xi$ (Cf. \cite{stoer:80}, pp.
84-88, \cite{gen:69}). He advised to avoid this technique,
particularly for low frequencies $\xi$ (e.g. for $\xi=0$).
However, we prove that under natural assumptions Goertzel's
algorithm is numerically stable in a sense (2), as an algebraic
polynomial evaluation algorithm. These results extend the results
obtained in \cite{gen:69}, \cite{smok:81} for real coefficients
$a_n$. Here we consider more general case of complex coefficients
$a_n$.

In the exact arithmetic we have for the quantities computed by
Goertzel's algorithm (Algorithm 2)
\begin{equation}
b_n= \sum_{k=n}^{N} { a_k \,{\vert z \vert}^{k-n} \,\,
U_{k-n}(t)},\,\,\, n=1,2, \ldots, N,
\end{equation}
\begin{equation}
u= \sum_{k=0}^{N} { a_k {\vert z \vert}^{k} \,\, T_{k}(t)}, \,\,
\,\, v= y\,\, \sum_{k=1}^{N} { a_k {\vert z \vert}^{k-1} \,\,
U_{k-1}(t)},
\end{equation}
\begin{equation}
t=\frac{x}{\mid z \mid}, \,\,\, t \in [-1,1]
\end{equation}
\noindent and $T_k(t)$ and  $U_k(t)$ are the Chebyshev polynomials
of the first kind and the second kind, respectively. They satisfy
the recurrence relations (Cf. \cite{szego:59})
\[
T_k(t)=2 t T_{k-1}(t)-T_{k-2}(t),\,\, \, U_k(t)=2 t
U_{k-1}(t)-U_{k-2}(t),\,\,\; k=2, \ldots
\]
\noindent with $T_0(t)= U_0(t)= 1$ and $T_1(t)=t,\,\, U_1(t)= 2
t$.

Moreover,
\[
T_k(t)= t U_{k-1}(t)- U_{k-2}(t), \,\,\; k=2,3, \ldots
\]
We remind very well known inequalities for $t \in [-1,1]$:
\begin{equation}
\vert T_k(t) \vert \leq 1, \,\, \vert U_{k}(t) \vert \leq k+1,
\,\,\; k=0,1, \ldots
\end{equation}

It's well known \cite{szego:59} that for $\vert t \vert < 1$
\[
T_k(t)=\cos k\theta, \hskip.3cm U_k(t)=\frac{\sin
(k+1)\theta}{\sin \theta},
\]
where
\[
t=\cos \theta, \hskip.3cm z=\vert z \vert e^{i\theta}, \hskip.3cm
\theta \in (0, \pi).
\]

\noindent Notice that
\begin{equation}
\frac{z^k}{\vert z \vert^k}=T_k(t)+i\frac{y}{\vert y
\vert}U_{k-1}(t) \hskip.3cm \mathrm{for}  \hskip.3cm
k=0,1,\dots,N.
\end{equation}

Now we analyze numerical behaviour of Goertzel's algorithm in
floating-point arithmetic. 

Let $\tilde w(z)= \tilde{u} + i \tilde{v}$, $\tilde{b}_n$,
$\tilde{p}, \tilde{q}$ denote the quantities computed numerically
in cfl (see section 3.2). We have
\[
\tilde{p}=p, \,\, \tilde{q}= q(1+\gamma), \,\,\, \vert \gamma
\vert\leq 2 \eps + \mathcal O(\eps^{\,2}).
\]
\noindent Therefore, $\tilde{b}_{N+1}=0,\,\, \tilde{b}_N=a_{N}$
and for $n=N-1,\ldots, 1$ we get
\begin{equation}
\tilde{b}_n= (a_n + \eta_n) + p\,\tilde{b}_{n+1} + q
\,\tilde{b}_{n+2},
\end{equation}
\begin{equation}
\tilde{u} = (a_0 + \eta_0) + x \,\tilde b_1 + q \,\tilde b_2
\end{equation}
\noindent where for $n=0,1,\dots,N$
\begin{equation}
\vert \eta_n \vert \leq K\,\eps\,(\vert a_n \vert + \vert z \vert
\vert\tilde{b}_{n+1}\vert +{\vert z \vert}^2
\vert\tilde{b}_{n+2}\vert), \hskip.3cm K=5+\mathcal O (\eps).
\end{equation}
\noindent The constant $5$ in (24) is overestimated, but this way
error analysis is simpler and the essential result is the same.

Further, we get
\begin{equation}
\tilde w(z)=(\tilde u+i\,\tilde v)\,(1+\delta_2), \hskip.3cm
\tilde{v}= y \,\tilde b_1 (1+ \delta_1), \,\,\; \vert \delta_1
\vert,\,\vert \delta_2 \vert \leq \eps.
\end{equation}

\noindent From this it follows that
\begin{equation}
\tilde b_n= \sum_{k=n}^{N} {( a_k + \eta_k)  {\vert z \vert}^{k-n}
\,\, U_{k-n}(t)},\,\,\, n=1,2, \ldots, N
\end{equation}
\noindent and
\begin{equation}
\tilde u= \sum_{k=0}^{N} { (a_k + \eta_k)  {\vert z \vert}^{k}
\,\, T_{k}(t)}.
\end{equation}


\noindent From (17), (18) and (20) it follows that
\begin{equation}
\vert u \vert \leq g_0, \hskip.2cm \vert b_n \vert \leq
(N-n+1)g_n, \hskip.3cm n=1, \dots, N,
\end{equation}
\noindent where
\begin{equation}
g_n=\sum_{k=n}^{N}{ \,\vert a_k \vert \,\vert z \vert ^{k-n}},
\hskip.3cm n=0,1, \ldots, N.
\end{equation}

We want to estimate the absolute error $\vert \tilde{w}(z)-w(z)
\vert$. Let's write $\tilde b_n$ as $\tilde b_n=b_n+e_n$, where
from (20), (26) and (27)
\begin{equation}
\vert e_n \vert \leq (N-n+1) \sum_{k=n}^{N}{ \,\vert \eta_k \vert
\,\vert z \vert ^{k}}.
\end{equation}

\noindent The formulae (24), (26)-(29) yield
\[
\vert \eta_k \vert \leq K\,\eps\,(\vert a_k \vert + (N-k)\vert z
\vert g_{k+1} + (N-k-1){\vert z \vert}^2 g_{k+2})+\mathcal O
(\eps^{\,2}).
\]

\noindent Thus
\begin{equation}
\vert \eta_k \vert \leq 2K\,\eps\,(N-k)g_{k} + \mathcal O
(\eps^{\,2}), \hskip.3cm \mathrm{for}\hskip.3cm k=0,1,\dots,N.
\end{equation}

\noindent Now write analogously $\tilde u = u + e_0$, where $\vert
e_0 \vert \leq \sum_{k=0}^{N}{ \,\vert \eta_k \vert \,\vert z
\vert ^{k}}$.

\noindent It's easy to check that $\sum_{k=0}^{N}{ \,g_k \,\vert z
\vert ^{k}} \leq (N+1)\,g_0$. From this and (31) we get
\begin{equation}
\sum_{k=0}^{N}{ \,\vert \eta_k \vert \,\vert z \vert ^{k}} \leq 2K
(N+1)N \eps g_0 + \mathcal O (\eps^{\,2}).
\end{equation}

\noindent Now let's rewrite (25) as
\begin{equation}
\tilde w(z)=(\tilde u+i\,y \,\tilde b_1) +\xi.
\end{equation}

\noindent It's easy to verify that
\begin{equation}
\vert \xi \vert \leq (2N+1)\eps g_0 + \mathcal O (\eps^{\,2}).
\end{equation}

\noindent Further from (21), (26) and (27) we get
\[
\tilde u + i\,y \,\tilde b_1 = \sum_{k=0}^{N}{\, (a_k
+\eta_k)\,z^k}.
\]

\noindent This and (33) yield
\[
\vert \tilde{w}(z)-w(z) \vert \leq \sum_{k=n}^{N}{ \,\vert \eta_k
\vert \,\vert z \vert ^{k}} + \vert \xi \vert.
\]

\noindent Combining this with (32) and (34) we get the inequality
\[
\vert \tilde{w}(z)-w(z) \vert \leq 2K (N+1)^2 \eps g_0 + \mathcal
O (\eps^{\,2}),
\]

\noindent which can be reformulated in the following

\begin{theor}
Assume that $cfl(z)= z$ and $cfl(a_n) = a_n$ for $n=0, \ldots, N$.
Let
\begin{equation}
2K (N+1)^2 \,\eps \leq 0.1,
\end{equation}
where $K$ is defined in {\rm(24)}.

\noindent Then Goertzel's algorithm for computing
$w(z)=\sum_{n=0}^{N}{a_n \,z^n}$ is componentwise backward stable,
i.e.
\begin{equation}
\tilde w(z)= \sum\limits_{n=0}^N\,\,
a_n(1+\Delta_n)\,z^n,\hskip4ex \vert \Delta_n\vert \leq A_N\,\eps
+ \mathcal O(\eps^{\,2}),
\end{equation}
where
\begin{equation}
A_N=2K (N+1)^2.\hskip1ex\domek
\end{equation}
\end{theor}

Notice that $A_N\approx 10 N^2$. Numerical tests in section 4
confirm that the constant $N^2$ is realistic.

\baselineskip=0.9\normalbaselineskip

\end{document}